\documentclass[10pt,leqno]{amsart}

\usepackage{amssymb,amsmath,amsthm}
\usepackage{bbm}
\usepackage[latin1]{inputenc}
\usepackage[english]{babel}
\usepackage{mathrsfs}
\usepackage[margin=3cm]{geometry}

\usepackage{graphicx}
\usepackage{stackrel}
\usepackage[toc,page]{appendix}
\usepackage{float} \usepackage[colorlinks, citecolor=blue]{hyperref}
\usepackage{enumerate, enumitem}

\usepackage{tikz}
\usepackage{mathdots}
\usepackage{yhmath}
\usepackage{cancel}
\usepackage{color}
\usepackage{siunitx}
\usepackage{array}
\usepackage{multirow}
\usepackage{textcomp,gensymb}
\usepackage{tabularx}
\usepackage{extarrows}
\usepackage{booktabs}
\usetikzlibrary{fadings}
\usetikzlibrary{patterns}
\usetikzlibrary{shadows.blur}
\usetikzlibrary{shapes}

\newcommand{\R}{\mathbb R}

\newcommand{\al}{\alpha}

\newcommand{\ve}{\varepsilon}

\newcommand{\sset}{\subseteq} 

\theoremstyle{plain}
\newtheorem{thm}{Theorem}[section]
\newtheorem{cor}[thm]{Corollary}
\newtheorem{lem}[thm]{Lemma}
\newtheorem{prop}[thm]{Proposition}

\theoremstyle{definition}

\newtheorem{rem}[thm]{Remark}

\theoremstyle{remark}

\definecolor{ceruleanblue}{rgb}{0.16, 0.32, 0.75}
\definecolor{cocoabrown}{rgb}{0.82, 0.41, 0.12}
\definecolor{crimsonglory}{rgb}{0.75, 0.0, 0.2}
\definecolor{aqua}{rgb}{0, 0.51, 0.5}

\begin{document}
	
	\title{Mountain pass frozen planet orbits in the helium atom model}
	\author{Stefano Baranzini,  Gian Marco Canneori and Susanna Terracini}
	\address{Dipartimento di Matematica ``G. Peano''
		\newline\indent
		Universit\`a degli Studi di Torino
		\newline\indent
		Via Carlo Alberto 10, 10123 Torino, Italy\\}
	\email{stefano.baranzini@unito.it}
	\email{gianmarco.canneori@unito.it}
	\email{susanna.terracini@unito.it}
	
	\date{\today}

	\date{\today}
	
	\begin{abstract}
	We seek frozen planet orbits for the helium atom through an application of the Mountain Pass Lemma to the Lagrangian action functional. Our method applies to a wide class of gravitational-like interaction potentials thus generalising the results in \cite{Cieliebak_variational}. We also let the charge of the two electrons tend to zero and perform the asymptotic analysis to prove convergence to a limit trajectory having  a collision-reflection singularity between the electrons.
	\end{abstract}
	
	\maketitle
	\begingroup
	\hypersetup{linkcolor=black}
	\tableofcontents
	\endgroup

	\section{Introduction}
 
    Frozen planet orbits are motions of a $1$-dimensional helium atom model, where the nucleus is fixed and both electrons move on the same side of a line.  These motions are periodic, in the sense that one electron keeps bouncing against the nucleus, whereas the second one slowly oscillates far from it (the  ``frozen planet''), running along a brake trajectory. The existence of these orbits has been in the focus of the recent mathematical literature (see \cite{Cieliebak_Langmuir,Cieliebak_Non_deg,Cieliebak_variational, Lei_Zhao_shooting}), also for its relevance in the semiclassical analysis of the helium atom model (cf \cite{MR1488438}). The quoted papers deal with the existence of periodic solutions of the following equations:
    \begin{equation}
    	\label{eq:ODE_helium_zhao}
    	\begin{cases}
    		\ddot{q}_1 = -\frac{2}{q_1^2}-\frac{1}{(q_2-q_1)^2}\\
    		\ddot{q}_2 = -\frac{2}{q_2^2}+\frac{1}{(q_2-q_1)^2}
    	\end{cases}.
    \end{equation}
    This system models the behaviour of two collinear electrons $(q_1,q_2)$ in the helium atom, assuming that the nucleus is fixed at the origin. Each of the two particles is subjected to an attractive force $-\frac{2}{q_i^2}$ towards the nucleus, and to a repulsive force $\pm \frac{1}{\left(q_2-q_1\right)^2}$, which pushes the electrons apart. While \cite{Lei_Zhao_shooting} adopts a shooting method, 
    \cite {Cieliebak_Langmuir,Cieliebak_Non_deg,Cieliebak_variational} propose a variational approach to search of such orbits. Intrigued  by the indirect and somewhat cumbersome approach pursued there, we wondered about the possibility of directly applying some Critical Point Theorem, paving the way for the treatment of a broader class of models.

	Indeed, the aim of this paper is to study periodic solutions of the following system of one dimensional second order non-linear equations:
	\begin{equation}
		\label{eq:ODE_helium}
		\begin{cases}
			\ddot{q}_1 = f'(\vert q_1 \vert)+g'(\vert q_2-q_1\vert )\\
			\ddot{q}_2 = f'(\vert q_2\vert )-g'(\vert q_2-q_1\vert )
		\end{cases},
	\end{equation}
    where $f$ and $g$ are real functions both having a singularity at the origin. Here, $f$ represents the attraction force to the nucleus and $g$ the repulsive interaction between the particles.

    In our paper we will deal with a much larger class of potentials than the one in \cite{Cieliebak_Langmuir,Cieliebak_Non_deg,Cieliebak_variational, Lei_Zhao_shooting} and we will make the following natural assumptions on the  functions $f,g\colon\R^+\to\R^+$ involved in \eqref{eq:ODE_helium}. We assume that $f,g \in C^2(\mathbb{R}^+)$ and:
    \begin{align}
    	\label{hyp:decaying}
    	& f(s),f'(s),g(s),g'(s)\to 0,\quad\text{as}\ s\to+\infty \\
    	\label{hyp:monotonicity_convexity}
    	& f'(s),g'(s)\le 0\ \text{and}\ g''(s),f''(s)\ge 0\; \forall s>0\\
    	\label{hyp:homogeneity}
    	& \exists \al \in (0,2) \ \text{such that}\ sf'(s)+\al f(s)\ge 0,\ sg'(s)+\al g(s)\leq 0\;, \forall s>0\\
		\label{hyp:mountain_pass}
    	& \exists \bar s>0 \ \text{such that} \; 0<g(\bar s)<f(\bar s).
    \end{align}
   
   Note that potentials  $f(s)=a/s^\alpha,\, g(s)=b/s^\beta$ always fulfil these assumptions provided $\alpha\in (0,2)$, $\beta\geq\alpha$ and, whenever $\alpha=\beta$, $b<a$. In particular, also system \eqref{eq:ODE_helium_zhao} does. The class of admissible potentials, however, is clearly much broader.
   
   Assumption \eqref{hyp:homogeneity} is a classical homogeneity condition already present in literature (see e.g. \cite{ACZ1987,ACZ1990}). Condition \eqref{hyp:mountain_pass} ensures that, independently on the initial position of the electrons, the attractive singularity prevails on the repulsive one, when the outer electron $q_2$ is far enough. 
   
  Notice that the first equation of \eqref{eq:ODE_helium} may be singular at the origin, even though, as $\alpha\in(0,2)$ with a weak force.  Therefore, as $q_1$ has one degree of freedom,  \emph{any bounded solution} has a collision with the origin. Thus, we will construct \emph{generalized} periodic solutions, for which we allow for collisions between the first electron and the nucleus. 
    For a given $T>0$, we seek solutions $(q_1,q_2)$ of \eqref{eq:ODE_helium} satisfying $q_1(t)<q_2(t)$ for all $t\in[0,T]$ and such that:
    \begin{equation}
    	\label{eq:initial_conditions}
    	\begin{cases}
    	   \dot{q}_1(0) =0 =\dot{q}_2(0)\\
    		q_1(T) =0 =\dot{q}_2(T)
    	\end{cases}.
    \end{equation}
    Thus, by reflecting the trajectory of $q_1$ after its collision with the origin, we obtain a generalized periodic solution of period $2T$, with a unique collision at $t=T$. Such a solution will be called a \emph{frozen planet orbit of period $2T$}. Our main results are the following:
    
      \begin{thm}
    	\label{thm:main_frozen}
    	For any functions $f,g$ satisfying assumptions \eqref{hyp:decaying}-\eqref{hyp:mountain_pass} and for any $T>0$ there exists a \emph{frozen planet orbit}, i.e., a generalized periodic solution of \eqref{eq:ODE_helium} with period $2T$, which satisfies the boundary conditions given by \eqref{eq:initial_conditions} and has no other singularities. 
    \end{thm}
    
    \begin{rem} Of course, this Theorem applies also to \eqref{eq:ODE_helium_zhao}. See also Corollary \ref{cor:energy} further ahead for the fixed energy problem.
    \end{rem}
        
    In the last section of the paper we investigate the behaviour of \eqref{eq:ODE_helium} in the case in which the repulsive singularity is damped. To model this situation, we introduce a parameter $\mu \in [0,1]$ and consider the following singularly perturbed system:
    	\begin{equation}
    		\label{eq:ODE_mu_intro}
    	\begin{cases}
    		\ddot{q}_1 = f'(\vert q_1 \vert)+\mu g'(\vert q_2-q_1\vert )\\
    		\ddot{q}_2 = f'(\vert q_2\vert )-\mu g'(\vert q_2-q_1\vert )
    	\end{cases}.
    \end{equation}
    When $\mu=0$ we see that \eqref{eq:ODE_mu_intro} decouples into two Kepler-like equations. We shall prove  the following result:
    \begin{thm}
    	\label{thm:main_brake}
    	As $\mu\to0$ the $\mu-$frozen planet orbits constructed in Theorem \ref{thm:main_frozen} converge uniformly on compacts subsets of $(0,T)$ to a function $q = (q_1,q_2).$ Each of the components of $q_i$ is a segment of the brake orbit of period $2T$ of the $f-$Kepler problem:
    	\begin{equation*}
    		\ddot{x} = f'(x).
    	\end{equation*}
        Moreover, $q_1(T) =0$ and $\dot{q}_1(0) =- \dot{q}_2(0)>0$. So, the limit trajectory has  a collision-reflection at $t=0$.
    \end{thm}
   
   \begin{rem}
   	Theorem \ref{thm:main_brake} can be framed in the context of \emph{degenerate billiards} introduced in \cite{bolotin_degenerate_billiards} (see also \cite{zbMATH06723100} and references therein). 	   	In our case, the singularity set is the positive diagonal  $\left\lbrace (q_1,q_2)\in \mathbb{R}^2:\right.$  $\left. q_1,q_2\geq 0,\ q_1=q_2\right\rbrace$ and the billiard trajectory is the periodic trajectory obtained by \emph{regularizing} the collision at $q_1=0$ of the limit curve $q$ obtained in Theorem \ref{thm:main_brake}. A possible application to \eqref{eq:ODE_mu_intro} of the results in \cite{bolotin_degenerate_billiards}, in order to prove existence of periodic solutions for small values of $\mu$, would involve regularising the singularity of $f$ and verifying an appropriate non-degeneracy condition on the limiting profile.	Then, trajectories of \eqref{eq:ODE_mu_intro} could be obtained as perturbations of trajectories of two decoupled Kepler problems. Consequently, according with the technique proposed there, it could be possible to construct branches of solutions of \eqref{eq:ODE_mu_intro} emanating from billiard trajectories reflecting on the singularity set.  However, although we believe it possible, in this paper we will not use the perturbative approach to prove the existence of  solutions for $\mu$ small, since we do develop a global variational approach, valid for all values of $\mu$ compatible with our assumptions.  Indeed, thanks to our global approach, we show that this branch can be continued for values of $\mu$ belonging to the whole interval $(0,1]$.    
   \end{rem}

   As said, the approach we follow in this paper is variational. We wish to characterize solutions of \eqref{eq:ODE_helium} as critical points of the Lagrangian action functional $\mathcal{A}$ defined on  the open set $\mathcal{U} = \{(q_1,q_2)\in H^1([0,T],\mathbb{R}^2): q_1\ne q_2 \}$ of the Hilbert space $H^1([0,T],\mathbb{R}^2)$. Indeed, non-collision solutions of \eqref{eq:ODE_helium} correspond to critical points of
	\begin{equation}
		\label{eq:def_true_action}
		\mathcal{A}(q_1,q_2) = \int_0^T \frac12  \left(\vert\dot{q}_1\vert^2+ \vert\dot{q}_2\vert^2\right)+ f(\vert q_1 \vert)+f(\vert q_2\vert)-g(\vert q_2-q_1\vert).
	\end{equation}
	Unfortunately, we must cope with the presence of singularities both on $f$ and $g$, so the notion of critical point must be suitably generalized. Collisions of the first electron with the nucleus, in particular, can not be avoided since, as already mentioned, solutions of \eqref{eq:ODE_helium} with initial conditions \eqref{eq:initial_conditions} are  expected to have collisions. In order to deal with this problem, we introduce a family of modified functionals $\mathcal{A}_{\ve_1,\ve_2}$ depending on two small parameters $\ve_1,\ve_2>0$. They are defined as follows:
    \begin{equation}
    	\label{eq:def_smoothed_action}
    	\mathcal{A}_{\ve_1,\ve_2}(q_1,q_2) = \int_0^T \frac12 \left(\vert\dot{q}_1\vert^2+		\vert\dot{q}_2\vert^2\right)+ f_{\ve_1}( q_1)+f_{\ve_1}( q_2)-g_{\ve_2}(\vert q_1-q_2 \vert),
    \end{equation}
    where $f_{\ve_1}$ is a smooth function which approximates the singularity at the origin, while  $g_{\ve_2}$ is a penalization of $g$ with a \emph{strong force} term concentrated in an $\ve_2$-neighbourhood of collisions (precise definitions are given in Section \ref{sec;approximating_functions}, and their main properties are listed in Lemma \ref{lemma:properties_f_g}). 

    For any $\ve_1,\ve_2$ small enough, we will seek critical points of $\mathcal{A}_{\ve_1,\ve_2}$ in the set
    \[
    \mathcal{D}: = \{(q_1,q_2)\in \mathcal{U}:q_1(0)=0,\, q_1<q_2\}.
    \] 
    From a variational perspective, such critical points are always saddles and existence is proved via a slight variant of the \emph{Mountain Pass Lemma}   (see Lemma \ref{lemma:abstract_mountain_pass}). Thanks to some suitable a priori estimates on the energy and on the $H^1$ norms of solutions, we prove that up to subsequence they converge uniformly to a solution of \eqref{eq:ODE_helium} satisfying \eqref{eq:initial_conditions}. Compared with the variational approach developed in \cite{Cieliebak_Non_deg,Cieliebak_variational}, ours has the advantage of avoiding the regularization argument and, being more direct, of giving more informations on the action level of the solutions, which will come in handy in the analysis of the asymptotics as $\mu\to 0$. Finally, we believe it is amenable to an extension to the multi-electron case, which we will deal with in a subsequent study.
    
    The structure of the paper is the following. We first introduce the smoothings $f_{\ve_1}$ and $g_{\ve_2}$ and then verify the hypotheses of Lemma \ref{lemma:abstract_mountain_pass}. In Section \ref{sec:ve_critical_points} we prove the existence of frozen planet orbits as critical points for $\mathcal{A}_{\ve_1,\ve_2}$. In Section \ref{sec:critical_points_A} we prove some a priori bounds on solutions and some of their qualitative properties. Finally, Section \ref{sec:perturbative} is devoted to the proof of Theorem \ref{thm:main_brake}.
   
	\section{Existence of critical points for $\mathcal{A}_{\ve_1,\ve_2}$}
	\label{sec:ve_critical_points}
	In this section we introduce, for $\ve_1,\ve_2>0$ small enough, the smoothed functional $\mathcal{A}_{\ve_1,\ve_2}$. After that, we establish the existence of \emph{frozen planet orbits} for each $\ve_1,\ve_2$ via a mountain pass lemma.
	
	\subsection{Construction and properties of $f_{\ve_1}$ and $g_{\ve_2}$}\label{sec;approximating_functions}
	
	We start by introducing the function $f_{\ve_1}$. They approximate the function $f$, smoothing the Keplerian singularity at $0$. Let $p_{\ve_1,f}$ be the line tangent to $f$ at $s = \ve_1$, namely:
	\[
	   p_{\ve_1,f}(s) = f(\ve_1)+f'(\ve_1)(s-\ve_1).
	\] 
	Let $\tilde{p}_{\ve_1,f}$ be the degree two polynomial having the same value and derivative at $0$ as $p_{\ve_1,f}$ and a maximum in $-\ve_1$, which reads:
	\[
	 \tilde{p}_{\ve_1,f}(s) =   f(\ve_1)-f'(\ve_1)\ve_1+f'(\ve_1)s+\frac{f'(\ve_1)}{2\ve_1}s^2.
	\] 
	Define also $f_{\ve_1}$ as follows:
	\begin{equation}
		\label{eq:def_f_1}
		f_{\ve_1}(s) = 
		\begin{cases}
			\begin{aligned}
			&f(s) &\text{ if } s\ge \ve_1 \\
			&p_{\ve_1,f}(s) &\text{ if } s \in [0,\ve_1] \\
			&\tilde{p}_{\ve_1,f}(s) &\text{ if } s \in [-\ve_1,0] \\
			&\tilde{p}_{\ve_1,f}(-\ve_1) &\text{ if } s \le -\ve_1
			\end{aligned}
		\end{cases}.
	\end{equation}
	Moreover, consider a smooth monotone decreasing function $\psi_{\ve_2}\colon\R\to\R$ such that:
	\begin{itemize} 
		\item $\psi_{\ve_2}\vert_{(-\infty,-\ve_2]} \equiv 1$ and  $\psi_{\ve_2}\vert_{[\ve_2,+\infty)}\equiv 0$;
		\item $\psi_{\ve_2}(0) = \frac{1}{2}$;
		\item $\psi_{\ve_2}$ is convex on $[0,\ve_2]$.
	\end{itemize}	
    Then, we define $g_{\ve_2}$ as:
	\begin{equation*}
		g_{\ve_2}(s) = g(s)+\frac{\psi_{\ve_2}(s)}{s^2},
	\end{equation*}
	where the second addendum is a \emph{strong force} perturbation of the original repulsive singularity. In the following result we collect the main properties of the functions $f_{\ve_1}$ and $g_{\ve_2}$.
	\begin{lem}
		\label{lemma:properties_f_g}
		For all $\ve_1,\ve_2>0$, the following properties hold:
		\begin{enumerate}[label=\roman*)]
			\item $f_{\ve_1} \in C^{1,1}(\mathbb{R})$;
			\item $f_{\ve_1}$ is monotone decreasing on $\mathbb{R}$ and convex on $[0,+\infty)$;
			\item  $s\,f'_{\ve_1}(s)+ \alpha f_{\ve_1}(s)\ge 0$ for all $s$;
			\item $g_{\ve_2}$ is monotone decreasing on $(0,+\infty)$ and convex;
			\item  $s\,g'_{\ve_2}(s)+ \alpha g_{\ve_2}(s)\le 0$ for all $s>0$.
		\end{enumerate}
		\begin{proof}
			Properties $i)$ and $ii)$ follow by construction. Let us prove point $iii)$. The inequality holds on $[\ve_1,\infty)$ since $f_{\ve_1}$ coincides with $f$ and on $(-\infty,-\ve_1]$ since $f_{\ve_1}$ is constant. When $s\in[0,\ve_1]$, $f_{\ve_1}$ coincides with $p_{\ve_1}$, so it is a consequence of assumptions \eqref{hyp:monotonicity_convexity}-\eqref{hyp:homogeneity} on $f$.  On the other hand, when $s\in[-\ve_1,0]$, a direct computation shows that:
			\begin{equation*}
				\left(\alpha \tilde{p}_{\ve_1,f}(s) +s 	\tilde{p}_{\ve_1,f}'(s)\right)' = (\alpha+1)	\tilde{p}_{\ve_1,f}'(s) + s \tilde{p}_{\ve_1,f}''(s)= (\al +1)f'(\ve_1)+(\al+2)\frac{f'(\ve_1)}{\ve_1}s,
			\end{equation*}
			and so the maximum of $\alpha \tilde{p}_{\ve_1,f}(s) +s 	\tilde{p}_{\ve_1,f}'(s)$ is achieved at $\bar{s}=-\frac{\al+1}{\al+2}\,\ve_1$, which is an internal point. Since $\alpha 	\tilde{p}_{\ve_1,f}(s) +s 	\tilde{p}_{\ve_1,f}'(s)$ is greater than the minimum of its values at $0$ and $-\ve_1$, which are both positive, we are done. 
			
			Point $iv)$ and $v)$ hold for $g$ by assumption so we have to check them for $\frac{\psi_{\ve_2}(s)}{s^2}$.
			A straightforward computation shows that for $s>0$:
			\[
			\begin{aligned}
				\left(\frac{\psi_{\ve_2}}{s^2} \right)' =& \frac{s\psi_{\ve_2}'-2\psi_{\ve_2}}{s^3}\le 0, \\
				\left(\frac{\psi_{\ve_2}}{s^2} \right)''=& \frac{s^2\psi_{\ve_2}''-4s\psi_{\ve_2}'+6\psi_{\ve_2}}{s^4}\ge 0,
			\end{aligned}
			\]
			proving $iv)$. Finally notice that:
			\[
			 s	\left(\frac{\psi_{\ve_2}}{s^2} \right)'+\alpha \frac{\psi_{\ve_2}}{s^2} = \frac{s(\psi_{\ve_2}'+\alpha \psi_{\ve_2})-2 \psi_{\ve_2}}{s^3} \le 0
			\]
			as soon as $s$ is small enough, since $2 \psi_{\ve_2}(0)=1$ and $\psi'_{\ve_2}$ is decreasing.
		\end{proof}
	\end{lem}

	\subsection{Mountain pass geometry}
	We will consider now the problem of finding critical points of $\mathcal{A}_{\ve_1,\ve_2}(q_1,q_2)$ on the set $\mathcal{D} = \{(q_1,q_2):\,q_1(0)=0,\ q_1<q_2\}$. The basic tool we will use is a mountain pass type theorem. In this section we show how to set up a mountain-pass type geometry. 
	
	For every $\ve_1>0$, denote by $q_{\ve_1}$ the brake orbit of the smoothed Kepler problem having period $T$. It coincides with the minimiser of:
	\begin{equation*}
		\mathcal{F}_{\ve_1} (q) = \int_0^T \frac12 \vert \dot{q }\vert^2+f_{\ve_1}(q)
	\end{equation*}
     on $\{q \in H^1([0,T]): q(0) =0\}$. Let $a_{\ve_1}$ be the value of $\mathcal{F}_{\ve_1}(q_{\ve_1})$ (see Lemma \ref{lem:brake} for details). Fix some constants $c,C\in \mathbb{R}$ such that $C>c>q_{\ve_1}(T_{\ve_1})$ and let us define the following family of continuous paths $\gamma :[0,1] \to \mathcal{D}$
	\begin{equation}
		\label{eq:min_max_classes}
		\Gamma_{c}^C(\ve_1)=\left\lbrace\gamma\colon[0,1] \to \mathcal{D}\text{ such that }\gamma (0) = (q_{\ve_1}, c) \text{ and }\gamma(1)= (q_{\ve_1}, C) \right\rbrace.
	\end{equation}
	To ease the notation, we will write $\gamma_s$ in place of $\gamma(s)$. Let us prove the following 
	\begin{lem}
		\label{lem:mountain_pass_geom}
		For $\ve_1,\ve_2$ small enough, there exist $\delta_1 >\delta_2>0$ and $c, C\in(q_{\ve_1}(T),+\infty)$ such that, for any path $\gamma \in \Gamma_c^C(\ve_1)$, we have
        \[
           \max\{\mathcal{A}_{\ve_1,\ve_2}(\gamma_0),\mathcal{A}_{\ve_1,\ve_2}(\gamma_1)\}\le a_{\ve_1}+\delta_2,
		  \quad
		   \max_{s \in [0,1]} \mathcal{A}_{\ve_1,\ve_2}(\gamma_s) > a_{\ve_1}+\delta_1.
        \]
		\begin{proof}
			The first step is to identify the value of $\mathcal{A}_{\ve_1,\ve_2}(\gamma_0)$ and $\mathcal{A}_{\ve_1,\ve_2}(\gamma_1)$. Recall that $f_{\ve_1}(s)=f(s)$ whenever $s>\ve_1$. For $c>\ve_1$, we have that
		   \[
			\begin{aligned}
				\mathcal{A}_{\ve_1,\ve_2}(q_{\ve_1},c) &= \int_0^T\frac12 \dot{q}_{\ve_1}^2 + f_{\ve_1}(q_{\ve_1})+f_{\ve_1}(c)-g_{\ve_2}(c-q_{\ve_1}) \\ &=\mathcal{F}_{\ve_1}(q_{\ve_1})+\int_0^Tf(c)-g_{\ve_2}(c-q_{\ve_1}).
			\end{aligned}
			\]
			Since $q_{\ve_1}$ solves the equation $\ddot{q}_{\ve_1} = f_{\ve_1}'(q_{\ve_1})$ on $[0,T]$, we can compute its Taylor expansion at $t = T$. Recalling that $\dot{q}_{\ve_1}(T)=0$, for $\ve_1$ sufficiently small this yields:
			\begin{equation}\label{eq:taylor_exp}
				q_{\ve_1}(t) = q_{\ve_1}(T) + f'(q_{\ve_1}(T)) (T-t)^2 + O((T-t)^3).
			\end{equation}
			Thus, $c-q_{\ve_1}(t) = c-q_{\ve_1}(T)+ f'(q_{\ve_1})(T-t)^2$. By definition of $g_{\ve_2}$, $g_{\ve_2}((T-t)^2)$ is not in $L^1$. It follows that			
			\begin{equation*}
				\liminf_{c\to q_{\ve_1}(T)} \int_0^Tg_{\ve_2}(c-q_{\ve_1}(t))dt \ge \int_0^Tg_{\ve_2}(q_{\ve_1}(T)-q_{\ve_1}(t))dt = +\infty,
			\end{equation*}
			and so:
			\[
			\lim_{c\to q_{\ve_1}(T)}\mathcal{A}_{\ve_1,\ve_2}(\gamma_0) = -\infty.
			\]
			Concerning the other endpoint, for $C$ large enough we have:
			\begin{equation*}
				\mathcal{A}_{\ve_1,\ve_2}(\gamma_1)=\mathcal{A}_{\ve_1,\ve_2}(q_{\ve_1},C)  = a_{\ve_1}+\int_0^Tf(C)-g(C-q_{\ve_1}).
			\end{equation*}
			Since both $f$ and $g$ go to zero at infinity thanks to \eqref{hyp:decaying}, $\mathcal{A}_{\ve_1,\ve_2}(\gamma_1)$ is arbitrarily close to $a_{\ve_1}.$ This concludes the proof of the first claim of our statement.
			
			Now we have to show that the maximum of $\mathcal{A}_{\ve_1,\ve_2}$ along any path is greater than $a_{\ve_1}+\delta_1$ for some $\delta_1$. For $d>0$, let us consider the set of paths $(q_1,q_2)\in \mathcal{D}$ at distance at least $d$:
			\begin{equation*}
				\mathcal{V}_d =\left\lbrace(q_1,q_2)\in\mathcal{D}:\min_{t\in[0,T]}\left(q_2(t)-q_1(t)\right) = d\right\rbrace.
			\end{equation*}
			Since any path $\gamma\in\Gamma_{c}^C(\ve_1)$ joins $(q_{\ve_1},c)$ to $(q_{\ve_1},C)$, by continuity, it has to cross $\mathcal{V}_d$  as soon as $d< C-q_{\ve_1}(T)$. The next step is to estimate the quantity:
			\begin{equation*}
				\min_{q \in \mathcal{V}_d} \mathcal{A}_{\ve_1,\ve_2}(q).
			\end{equation*}
			Without loss of generality, let us further assume that $d>\max\{\ve_1,\ve_2\}$. Since $\min_{t \in [0,T]} (q_2(t)-q_1(t)) = d$, we have  $f_{\ve_1}(q_2)=f(q_2)$ and $g_{\ve_2}(q_2-q_1)=g(q_2-q_1)$.  Let us denote by $q$ an element $q = (q_1,q_2) \in \mathcal{V}_d$. 
			
			First of all notice that, as $d$ grows, the first component $q_1$ of any minimiser must remain bounded. Indeed, $\mathcal{A}_{\ve_1,\ve_2}(q) \ge \frac12 \Vert \dot{q_1}\Vert^2_2	- T \,g(d)$ and, setting $D = d+q_{\ve_1}(T)$, we have:
			\[
			 \min_{q \in \mathcal{V}_d} \mathcal{A}_{\ve_1,\ve_2}(q)\le \mathcal{A}_{\ve_1,\ve_2}(q_{\ve_1},D)  = a_{\ve_1}+ \int_0^T f(D)-g(D-q_{\ve_1}(t)) \le a_{\ve_1}+T(f(D)-g(d)).		
			\]
			Let assume that the quantity $f(q_2)-g(q_2-q_1)\le0$ at some instant $t^*$. We claim that this implies that $q_1(t^*)$ must be large as $d$ goes to infinity.
			Indeed, thanks to assumption \eqref{hyp:homogeneity}-\eqref{hyp:mountain_pass}, we have the following two inequalities on $f$ and $g$ for $s\ge \bar{s}$:
			\[
			f(s) \ge f(\bar{s}) \left(\frac{\bar{s}}{s}\right)^\alpha \quad g(s) \le g(\bar{s}) \left(\frac{\bar{s}}{s}\right)^\alpha.
			\]
			As soon as $q_2-q_1$ and $q_2$ are greater than $\bar{s}$, this implies that:
			\[
			 	f(q_2)-g(q_2-q_1) \ge  f(\bar{s}) \left(\frac{\bar{s}}{q_2}\right)^\alpha-g(\bar{s}) \left(\frac{\bar{s}}{q_2-q_1}\right)^\alpha.
			\]
			The left-hand side is negative in $t^*$, so setting $\nu_{\alpha} = \left(\frac{f(\bar{s})}{g(\bar{s})}\right)^{\frac{1}{\alpha}}>1$ and evaluating at $t=t^*$ we obtain the following inequality:
		    \begin{equation*}
		    	q_1(t^*)\ge \frac{\nu_\alpha-1}{\nu_\alpha} q_2(t^*).
		    \end{equation*} 
	        Since $\min_{t \in[0,T]}q_2(t)\ge d$ and $q_1$ is bounded, $f(q_2)-g(q_2-q_1)$ is always positive. 		
			
		    Summarising, we have proved that for $d$ large enough, minimisers in $\mathcal{V}_d$ satisfy 	$\max_{t} \vert q_1(t)\vert\le q^*$ and  $f(q_2)-g(q_2-q_1)>0$. Let $q_d = (q_1^d,q_2^d) \in \mathcal{V}_d$ a minimiser. We have:
			\begin{equation*}
				\mathcal{A}_{\ve_1,\ve_2}(q_d) \ge \min_{q_1:q_1(0) =0} \mathcal{F}_{\ve_1}(q_1)+\int_0^T\frac12 \vert \dot{q}_2^d\vert+f(q_2^d)-g(q^d_2-q_1^d)\ge a_{\ve_1}+\int_0^Tf(q_2^d)-g(q^d_2-q_1^d).
			\end{equation*}
			On the other hand $\int_0^T f(q_2^d)-g(q^d_2-q_1^d)\ge T\min_{t\in [0,T]} (f(q_2^d)-g(q^d_2-q_1^d))>0$ proving the Lemma. 
		\end{proof}
	\end{lem}
    
    The mountain pass lemma requires also that the functional $\mathcal{A}$ is unbounded on the boundary $\partial\mathcal{D}$ (see Lemma \ref{lemma:abstract_mountain_pass}). Notice that
    \[
    	\partial\mathcal{D}=\{(q_1,q_2)\in H^1([0,T],\R^2):\ q_1(0)=0, \, q_1(t)=q_2(t),\ \text{for some}\ t\in[0,T]\}
    \]
    and, by construction, $g_{\ve_2}$ behaves like a strong force close to $\partial D$. The following Lemma is a straightforward modification of a classical argument needed to show that, in the strong force case, the action blows up at collisions (see e.g. \cite[Lemma 5.3]{AmbCotZel1993}).
    \begin{lem}
    	\label{lemma:explode_on_boundary}
    	If $q \in \partial \mathcal{D}$ then $\mathcal{A}_{\ve_1,\ve_2}(q)=-\infty$.
    	\begin{proof}
    		By construction, $g_{\ve_2}(s)\ge \psi_{\ve_2}(s)/s^2$. It follows that, for any $[t,s]\subseteq[0,T]$:
    		\begin{equation*}
    			\int_0^T g_{\ve_2}(q_2-q_1)\ge \int_t^s\frac{\psi_{\ve_2}(q_2-q_1)}{(q_2-q_1)^2}.
    		\end{equation*}
    		Let $q =(q_1,q_2) \in \partial\mathcal{D}$. Assume without loss of generality that $q_1(s) = q_2(s)$ and $q_1(w)<q_2(w)$ for all $w \in[t,s)$. Up to choosing a bigger $t$, we can assume that $\psi_{\ve_2}(q_2-q_1)\ge 1/4$.  We have:
    		\begin{align*}
    			\vert\log((q_2-q_1)(w))-	\log((q_2-q_1)(t))\vert &\le \int_t^w\frac{\vert\dot{q}_2-\dot{q}_1\vert}{q_2-q_1}\\
    			&\le\left(\int_t^w\frac{1}{(q_2-q_1)^2}\right)^2 (\Vert \dot q_1 \Vert^2+\Vert \dot q_2\Vert^2)^2\\
    			&\le 16 \left(\int_0^T g_{\ve_2}(q_2-q_1)\right)^2 \Vert \dot{q}\Vert^2_2.
    		\end{align*} 
    		Taking the limit as $w\to s$ we obtain that $g_{\ve_2}\notin L^1[0,T]$. Since $f_{\ve_1}$ is bounded we conclude that $\mathcal{A}_{\ve_1,\ve_2}(\partial \mathcal{D}) = -\infty.$ 
    	\end{proof}
    \end{lem}

    \subsection{Palais-Smale condition}
    It remains to show that the action functional $\mathcal{A}_{\ve_1,\ve_2}$ satisfies the Palais-Smale condition.    
    Let us define the candidate critical value of $\mathcal{A}_{\ve_1,\ve_2}$ as:
    \begin{equation}
    	\label{eq:min_max_value}
    	c_{\ve_1,\ve_2} = \min_{\gamma\in \Gamma_c^C(\ve_1)} \max_{s\in[0,1]} \mathcal{A}_{\ve_1,\ve_2}(\gamma_s)> a_{\ve_1}.
    \end{equation}
    We recall that $(u_n)\sset\mathcal{D}$ is a Palais-Smale sequence, (PS) in short, at level $c_{\ve_1,\ve_2}$ if $\mathcal{A}_{\ve_1,\ve_2}(u_n) \to c_{\ve_1,\ve_2}$ and $d_{u_n}\mathcal{A}_{\ve_1,\ve_2} \to 0$ in the $H^{-1}$ norm. 
    \begin{prop}
    	\label{prop:palais_smale_condition}
    	Any Palais-Smale sequence at level $c_{\ve_1,\ve_2}$ in $\mathcal{D}$ admits a strongly convergent subsequence.
    	\begin{proof}
    		First of all, let us show that any Palais-Smale sequence is bounded in $H^1$. To ease the notation, for $q=(q_1,q_2)\in\mathcal{D}$, define the total potential
    		\[
    		U(q)=f_{\ve_1}(q_1)+f_{\ve_1}(q_2)-g_{\ve_2}(q_2-q_1).
    		\]
    		For a (PS) sequence $(u_n)=(q_1^n,q_2^n)\sset\mathcal{D}$ and a test function $v$, we can compute
    		\begin{equation}
    			\label{eq:gradient_A_ve}
    		    \begin{aligned}
    			    d_{u_n}\mathcal{A}_{\ve_1,\ve_2}(v) = \int_0^T \langle \dot{u}_n,\dot{v}\rangle +\langle \nabla U(u_n),v\rangle,\\
    			    \nabla U(u_n) = \begin{pmatrix}
    				    f_{\ve_1}'(q^n_1)\\
    				    f_{\ve_2}'(q^n_2)
    			    \end{pmatrix}+\begin{pmatrix}
    				    g_{\ve_2}'(q^n_2-q^n_1)\\
    				    -g_{\ve_2}'(q^n_2-q^n_1)
    			    \end{pmatrix}.
    		    \end{aligned}
    		\end{equation}
    		In particular, choosing $v =u_n$ and using point $iii)$ and $v)$ of Lemma \ref{lemma:properties_f_g}, we have that:
    		\[
    		\begin{aligned}
    			\int_0^T\langle	\begin{pmatrix}
    				f_{\ve_1}'(q^n_1)\\
    				f_{\ve_2}'(q^n_2)
    			\end{pmatrix},\begin{pmatrix}
    				q_1^n \\q_2^n
    			\end{pmatrix} \rangle &= \int_0^T 	f_{\ve_1}'(q^n_1)q_1^n+
    			f_{\ve_2}'(q^n_2)q^n_2 \\
    			&\ge -\alpha\int_0^T f_{\ve_1}(q_1^n)+f_{\ve_1}(q_2^n)
    		\end{aligned}
    		\]
    		and
    		\[
    		\begin{aligned}
    			\int_0^T\langle	\begin{pmatrix}
    				g_{\ve_2}'(q_2^n-q_1^n)\\
    				g_{\ve_2}'(q^n_2-q_1^n)
    			\end{pmatrix},\begin{pmatrix}
    				q_1^n \\q_2^n
    			\end{pmatrix} \rangle &= -\int_0^Tg_{\ve_2}'(q_2^n-q_1^n) (q_2^n-q_1^n) \\
    			&\ge  \alpha \int_0^T g_{\ve_2}(q_2^n-q_1^n).
    		\end{aligned}
    		\]
    	   It follows that:
    		\begin{equation}
    			\label{eq:estimate1}
    			d_{u_n}\mathcal{A}_{\ve_1,\ve_2}(u_n)\ge \Vert \dot{u}_n \Vert_2^2- \alpha \int_0^T U(u_n)=\frac{2+\alpha}{2}\Vert\dot{u}_n\Vert_2^2- \alpha \mathcal{A}_{\ve_1,\ve_2}(u_n).
    		\end{equation}
    		Since $(u_n)$ is a (PS) sequence, we know that
    		\begin{equation*}
    			\mathcal{A}_{\ve_1,\ve_2}(u_n) =  \frac{1}{2}\Vert\dot{u}_n \Vert^2_2+\int_0^T U(u_n)\to c_{\ve_1,\ve_2}>0.
    		\end{equation*}
    		Assume by contradiction that $(\dot{u}_n)$ is unbounded. Combining with \eqref{eq:estimate1}, for $n$ large enough we obtain:
    		\begin{align*}
    			d_{u_n}\mathcal{A}_{\ve_1,\ve_2}(u_n) >0.
    		\end{align*}
    		Being a (PS) sequence implies also that $d_{u_n}\mathcal{A}_{\ve_1,\ve_2}(u_n)/\Vert u_n\Vert_{H^1}\to 0$ and so we deduce that
    		\begin{equation}
    			\label{eq:norm_proof_PS}
    			\frac{\Vert\dot{u}_n\Vert_2^2}{\Vert u_n \Vert_{H^1}} \to 0.
    		\end{equation}
    		Up to a subsequence, assume that $\Vert\dot{u}_n\Vert_2\to+\infty$. In order for \eqref{eq:norm_proof_PS} to hold, we must have that
    		\begin{equation}
    			\label{eq:position_velocity}
    			\frac{\Vert{u}_n\Vert_2}{\Vert\dot{u}_n\Vert_2^2}\to +\infty.
    		\end{equation}
    		Since $q_1^n(0)=0$, a Poincar\'e inequality holds for $q_1^n$ and so $\Vert q_1^n\Vert_2\le \sqrt{T}\Vert \dot{q}_1^n\Vert_2$.
    		
    		Note that a \emph{modified} Poincaré inequality holds for $q_2^n$ as well. Let us define the function $x_n(t) =q_2^n(T-t)$ which, for any $s,w\in[0,T]$, using Jensen inequality, satisfies:
    		\begin{equation*}
    			(x_n(s)-x_n(w))^2 =\left( -\int_w^s\dot{q}_2^n(T-t)\right)^2 \le  \int_0^T(\dot{q}_2^n(T-t))^2.
    		\end{equation*}
    		Integrating with respect to $s$, this implies that:
    		\begin{equation*}
    			\Vert q_2^n-q_2^n(w)\Vert_2^2 = \int_0^T \left (q_2^n(s)-q_2^n(w)\right)^2\,ds\le T \Vert \dot{q}_2^n\Vert^2_2
    		\end{equation*}
    		Thus, for any $w\in[0,T]$ and for a positive constant $C$ depending only on $T$, we have obtained
    		\[
    		\Vert q_2^n\Vert_2 \le C \left(\Vert \dot{q}_2^n\Vert_2+\vert q_2^n(w)\vert\right).
    		\] 
    		On the other hand, the unboundedness of $\Vert \dot{u}_n\Vert$ implies that there exists some $(w_n)\sset[0,T]$ such that $q_2^n(w_n)-q_1^n(w_n)\to0 $. Indeed, since both the value of the action and $f_{\ve_1}$ are bounded over the (PS) sequence, there exists a constant $C>0$ such that
    		\[
    			\left\lvert \frac{1}{2} \Vert \dot{u}_n \Vert^2_2-\int_0^T g_{\ve_2}(q_2^n-q_1^n)\right\rvert \le C.
    		\] 
    		Thus $\Vert \dot{u}_n \Vert^2_2$ explodes if and only if $q_1^n$ and $q_2^n$ get closer and closer. 
    	    Thus, for some constant $C_1>0$ we have
    		\[
    		\Vert q_1^n\Vert_2 +\Vert q_2^n\Vert_2 \le \sqrt{T} \Vert \dot{q}_1^n\Vert_2 + C\left(\Vert \dot{q}_2^n\Vert_2+\vert q_2^n(w)\vert\right)\le C_1\left(\Vert \dot{q}_1^n\Vert_2 + \Vert \dot{q}_2^n\Vert_2\right)
    		\] 
    		and so $\Vert u_n \Vert_2$ is controlled by $\Vert \dot{u}_n\Vert_2$, a contradiction for \eqref{eq:position_velocity}.
    		
    		So far, we have shown that $\Vert \dot{u}_n\Vert_2$ and $\Vert q_1^n \Vert_2$ are bounded. To prove the $H^1$ boundedness of (PS) sequences we have to show that $q_2^n$ is bounded in $L^2$ too. Let assume by contradiction that this is not the case. There exists a sequence $c_n\to +\infty$ such that $q_2^n\ge c_n$ and this implies that
    		\begin{equation}
    			\label{eq:estimate3}
    		\nabla U(u_n)-\begin{pmatrix}
    			f_{\ve_1}'(q_1^n) \\
    			0
    		\end{pmatrix}\to 0\quad\text{uniformly as}\ n\to+\infty.
    		\end{equation}
    	    Testing $d_{u_n}\mathcal{A}_{\ve_1,\ve_2}$ on $(0,q_2^n-q_2^n(0))$, one finds that $\Vert \dot{q}_2^n\Vert_2 \to 0$ and so $q_2^n-q_2^n(0)$ converges strongly to $0$.
    	    
    		Moreover, let $\mathcal{F}_{\ve_1}$ be the functional defined in \eqref{eq:smoothed_kepler}. Then, for any $v \in H^1$, $v(0)=0$ we have 
    		\[
    		d_{q_1^n}\mathcal{F}_{\ve_1}(v)=\langle \dot{q}_1^n,\dot{v}_1\rangle +\int_0^T f_{\ve_1}'(q_1^n)v_1, \text{ and } \Vert d_{q_1^n}\mathcal{F}_{\ve_1}\Vert\to0.
    		\]	
    		Thus, we obtained that $(q_1^n)$ is a (PS) sequence for $\mathcal{F}_{\ve_1}$. Thanks to Lemma \ref{lem:brake}, it converges to a brake orbit $q_{\ve_1}$ for a smoothed Kepler problem and so we conclude that 
    		\[
    		\mathcal{A}_{\ve_1,\ve_2}(q_1^n,q_2^n)\to a_{\ve_1}<c_{\ve_1,\ve_2},
    		\]
    		a contradiction.
    		
    		So far, we have proved that $u_n$ is bounded in $H^1$ and that
    		\begin{equation}
    			\label{eq:estimate4}
    		\min_{t\in [0,T]}\vert q_1^n(t)-q_2^n(t)\vert>d_\ve.
    		\end{equation}
    		Up to a subsequence, $u_n$ admits a weak limit $u$. In particular $u_n\to u$ uniformly and in $L^2$. 
    		Thanks to uniform convergence and the bound on the distance between $(q_1^n,q^n_2)$, dominated convergence implies that:
    		\begin{equation*}
    			\int_0^T\langle\nabla U(u_n),(u-u_n) \rangle \to 0
    		\end{equation*}
    		By hypothesis $d_{u_n}\mathcal{A}_{\ve_1,\ve_2}\to 0$, and so we obtain strong convergence since:
    		\begin{equation*}
    			d_{u_n}\mathcal{A}_{\ve_1,\ve_2}(u-u_n) = \langle \dot{u}_n,\dot{u}-\dot{u}_n \rangle+o(1) \to 0.
    		\end{equation*} 
    	\end{proof}
    \end{prop}
    
    \subsection{Existence of critical points}
    
    We are now in the position to prove this result:
    \begin{prop}
    	\label{prop:existence_scritical_points_smoothed}
    	For any $\ve_1,\ve_2>0$ small enough and any $T>0$, the functional $\mathcal{A}_{\ve_1,\ve_2}$ has a critical point at level $c_{\ve_1,\ve_2}$ in the set:
    	\begin{equation*}
    		\mathcal{D} = \{(q_1,q_2)\in H^1([0,T];\R^2):q_1<q_2,q_1(0)=0\}.
    	\end{equation*}
        In particular, there exists a collision-less solution of:
        \begin{equation}
       	\label{eq:motion_ve_critical_points}
       	\begin{cases}
       		\ddot{q}_1 = f'_{\ve_1}(q_1)+g'_{\ve_2}(q_2-q_1)\\
       		\ddot{q}_2 = f'_{\ve_1}(q_2) - g'_{\ve_2}(q_2-q_1)
       	\end{cases},
       \end{equation}
      satisfying $q_1(0) = \dot{q}_1(T) =0$ and $\dot{q}_2(0)=\dot{q}_2(T) =0$.
      \begin{proof}
      	This is an application of Lemma \ref{lemma:abstract_mountain_pass}. In Lemma \ref{lem:mountain_pass_geom} we proved that $c_{\ve_1,\ve_2}>a_{\ve_1}$, for a suitable choice of $C$ and $c$ in the definition of $\Gamma_{c}^C(\ve_1)$. Moreover, in Proposition \ref{prop:palais_smale_condition} we have shown that the (PS) condition holds at level $c_{\ve_1,\ve_2}.$ Finally, Lemma \ref{lemma:explode_on_boundary} verifies that $\mathcal{A}_{\ve_1,\ve_2}(\partial \mathcal{D}) = -\infty$.
      \end{proof}
    \end{prop}
    
    \section{Existence of critical points for $\mathcal{A}$}
    \label{sec:critical_points_A}
    In this section we exploit a limit argument in order to prove that our main problem \eqref{eq:ODE_helium} actually admits solutions. To this extent, the first step is to provide suitable a priori estimates on the $H^1$ norm and on the energy of each solution of \eqref{eq:motion_ve_critical_points}. 
    \subsection{A priori estimates and finer properties of solutions}
   
    To ease the notations we will denote solutions of \eqref{eq:motion_ve_critical_points} by $q^\ve = (q^\ve_1,q_2^\ve)$. In this first result, we detect some useful properties on the monotonicity of such solutions and their derivatives.
    \begin{lem}
    	\label{lemma:convexity_ve_solutions}
    	For every $\ve_1,\ve_2>0$, the following hold:
    	\begin{enumerate}[label =\roman*)]
    		\item $q_1^{\ve}$ and $q_2^{\ve}+q_1^\ve$ are concave. In particular they are positive for all $t\le T$ and admit a maximum at $T$;
    		\item $q^\ve_2-q^\ve_1$ is convex with a minimum at $t=T$;
    		\item $\vert\dot q^\ve_2\vert \leq \vert\dot q^\ve_1\vert$;
    		\item $q^\ve_1$ is monotone increasing and $q^\ve_2$ is monotone decreasing.
    	\end{enumerate} 
    	\begin{proof}
    		By the properties listed in Lemma \ref{lemma:properties_f_g} and  \eqref{eq:motion_ve_critical_points}, $q^\ve_1$ and $q^\ve_2+q^\ve_1$ are concave. Thus they are monotone with a maximum point at $t=T$ since $\dot{q}_1^\ve(T) = 0 =\dot{q}^\ve_2(T).$ Since $f_{\ve_1}$ is convex on $[0,\infty)$, $f'_{\ve_1}$ is increasing and $f'_{\ve_1}(q_2)-f'_{\ve_2}(q_1)\geq 0$ provided $q_2\geq q_1$. It follows that  $\ddot{q}_2^\ve-\ddot{q}_1^\ve\geq0$ and $q^\ve_2-q^\ve_1$ is convex. By the boundary conditions, $t=T$ is a critical point and thus a minimum.
    		
    		Assertion $iii)$ follows from the fundamental theorem of calculus. Indeed:
    		\begin{equation*}
    			\dot{q}^\ve_i(t) = -\int_t^T\ddot{q}^\ve_i(s)ds, \,\, i=1,2,
    		\end{equation*}
    		and, observing that  $\ddot{q}^\ve_1+\ddot{q}^\ve_2\le0$ and $\ddot{q}^\ve_2-\ddot{q}^\ve_1\ge0$, we conclude that $-\dot{q}^\ve_1\leq\dot{q}^\ve_2\leq\dot{q}_1^\ve$.
    		
    		The last assertion is proved as follows. We already know that $q_1^\ve$ is strictly increasing since it is strictly concave. Let assume that there is a point $t^*<T$ which is a strict minimum for $q_2^\ve$. In particular, $q^\ve_2(t)>q^\ve_2(t^*)$ for all $t>t^*$ small enough and $\ddot{q}_2^\ve(t^*)\ge0$. By Lemma \ref{lemma:properties_f_g}, $f'_{\ve_1}$ is monotone increasing, and so 
    		\[
    		f'_{\ve_1}(q_2^\ve(t))-f'_{\ve_1}(q_2^\ve(t^*))\ge0, \text{ if } q_2^\ve(t)>q_2^\ve(t^*).
    		\]
    		Similarly, since $q_2^\ve-q_1^\ve$ is decreasing (and this follows from $ii)$), we have that $-g_{\ve_2}'((q_2^\ve-q_1^\ve)(t))\ge-g_{\ve_2}'((q_2^\ve-q_1^\ve)(t^*))$ and so, plugging in \eqref{eq:motion_ve_critical_points}, we find that $\ddot{q}_2^\ve(t)\ge 0$ as long as $q_2^\ve(t)\ge q_2^\ve(t^*)$. Thus $q_2$ is convex on $[t^*,T]$ and there can be no stationary point at $t =T$. A contradiction.
    	\end{proof}
    \end{lem}
    The following proposition deals with the boundedness of solutions of \eqref{eq:motion_ve_critical_points} and their energies. 
    \begin{prop}
    	\label{prop:a_priori_estimates_energy}
    	Let $q^\ve=(q_1^\ve,q_2^\ve)$ be a critical point of $\mathcal{A}_{\ve_1,\ve_2}$ at level $c_{\ve_1,\ve_2}$ and let $h_{\ve}$ be the corresponding total energy value. Then, $q^\ve$ is uniformly bounded in $H^1$ and $h_{\ve}$ is uniformly bounded in $\mathbb{R}$, for $\ve_1,\ve_2$ sufficiently small. Moreover, we have:
    	\[
    	-\frac{c_{\ve_1,\ve_2}}{T}\le h_\ve \le\frac{(\alpha-2) \, c_{\ve_1,\ve_2}}{(2+\alpha)T}<0.
    	\]
    	\begin{proof}
    		Since $q^{\ve}$ is a critical point, $d_{q^\ve} \mathcal{A}_{\ve_1,\ve_2} =0$. In particular we have that:
    		\begin{equation*}
    			d_{q_{\ve}} \mathcal{A}_{\ve_1,\ve_2}(q^\ve) = \Vert \dot{q}^\ve \Vert_2^2+ \int_0^T \left(q_1^\ve f'_{\ve_1} (q_1^\ve) +q_2^\ve  f'_{\ve_1} (q_2^\ve)- (q_2^\ve-q_1^\ve)g'_{\ve_2}(q_2^\ve-q_1^\ve)\right)dt. 
    		\end{equation*}
    		By points $iii)$ and $v)$ of Lemma \ref{lemma:properties_f_g} we have that:
    		\[
    			0 = 	d_{q_{\ve}} \mathcal{A}_{\ve_1,\ve_2}(q^\ve) \ge \Vert \dot{q}^\ve \Vert_2^2-\alpha \int_0^T\left( f_{\ve_1}(q_1^\ve)+f_{\ve_1}(q_2^\ve)-g_{\ve_2}(q_2^\ve-q_1^\ve)\right)dt.
    		\]
    	    Rewriting the inequality as in the proof of Proposition \ref{prop:palais_smale_condition} (see in particular \eqref{eq:estimate1}), we obtain that:
    	    \begin{equation}\label{eq:critical_point_proof}
    	       \Vert \dot{q}^\ve\Vert_2^2 \le \frac{2\alpha \, c_{\ve_1,\ve_2}}{2+\alpha}.
    	    \end{equation}
    	    Clearly, $c_{\ve_1,\ve_2}$ is uniformly bounded in $\ve_1,\ve_2$ and so is $\Vert \dot{q}^\ve\Vert_2^2.$
    		
    	    Reasoning as in Proposition \ref{prop:palais_smale_condition}, $q^\ve_1$ is bounded in $H^1$ if and only if $\dot{q}_1^\ve$ is bounded in $L^2$. It remains to show that $q_2^\ve$ is bounded in $L^2$. Indeed, assume by contradiction that $q_2^\ve(T)\to \infty$. Testing $d_{q^\ve}\mathcal{A}_{\ve_1,\ve_2}$ on the variation $(0,q_2^\ve)$ we find that $\dot{q}_2^\ve$ goes to zero in $L^2$ and thus, as in Proposition \ref{prop:palais_smale_condition}, $c_{\ve_1,\ve_2}$ approaches the level of a brake orbit.  A contradiction.
    		
    	    Let us prove the bound on the energy. Integrating over $[0,T]$ we have that:
    	    \begin{equation*}
    	     	T h_\ve =\frac{1}{2}\Vert\dot{q}^\ve\Vert^2_2 -\int_0^T U(q^\ve) dt = \Vert\dot{q}^\ve\Vert^2_2-\mathcal{A}_{\ve_1,\ve_2}(q^\ve) = \Vert\dot{q}^\ve\Vert^2_2-c_{\ve_1,\ve_2}.
        	\end{equation*}
    		On the other hand, by equation \eqref{eq:critical_point_proof}, we have that:
    		\[
    		\frac{(\alpha-2) \, c_{\ve_1,\ve_2}}{2+\alpha}\ge \Vert \dot{q}^\ve \Vert_2^2- c_{\ve_1,\ve_2} = Th_\ve.
    		\]
    		Finally,  we have that $ h_\ve \ge -\frac{1}{T}c_{\ve_1,\ve_2}$.
    	\end{proof}
    \end{prop}
   We need to guarantee that we are not approaching a total collision, where both the electrons collapse into the nucleus. This is proven in the following:

    \begin{prop}[No total collision]
    	\label{prop:total_collision}
    	For any $\ve_1,\ve_2>0$ there exist constants $C_1,C_2>0$ not depending on $\ve_1,\ve_2$ such that:
    	\[
    	\Vert \dot{q}^\ve \Vert_2 \ge C_1,\quad q_1^\ve(T)\ge C_2.
    	\]
    	\begin{proof}
    	    Let us show that, if  $\dot{q}_1^\ve$ goes to zero in $L^2$, solutions $q^\ve$ of \eqref{eq:motion_ve_critical_points}  converge uniformly to zero. From point $iii)$ of Lemma \ref{lemma:convexity_ve_solutions}, we easily see that $\dot{q}_2^\ve$ to 0 in $L^2$ as well. Moreover, since $q_1^\ve(0)=0$, the Poincaré inequality shows that if $\dot{q}_1^\ve$ converges to zero in $L^2$, then $q_1^\ve(T)$ goes to zero as well.  Thus $q_2^\ve$ converges uniformly to a constant curve. The energies of solutions are uniformly bounded by Proposition \ref{prop:a_priori_estimates_energy}, and at $T$ read
    		\[
    		h_\ve = -f_{\ve_1}(q_1^\ve(T))-f_{\ve_1}(q_2^\ve(T)) + g_{\ve_2}(q_2^\ve(T)-q_1^\ve(T));
    		\]
    		therefore, $q_2^\ve(T)$ must converge to zero, otherwise $h_\ve\to-\infty$. So $q^\ve_1$ and $q^\ve_2$ converge uniformly to 0.  
    		
    		Having proved this claim, let us assume by contradiction that $q_1^\ve$ converges uniformly to 0. 
        	We have that:
        	\begin{equation*}
        	\begin{aligned}
    	    	\dot{q}_1^\ve(t) &= \int_t^T-f'(q_1^\ve)-g'_{\ve_2}(q_2^\ve-q_1^\ve) \\
    	    	q_1^\ve(T) &= -\int_0^T\int_0^T \chi_{[t,T]}(s)\left( f'(q_1^\ve)(s)+g'_{\ve_2}(q_2^\ve-q_1^\ve)(s)\right) ds dt
        	\end{aligned}
    	   \end{equation*}
        	where $\chi_{[t,T]}$ denotes the characteristic function of $[t,T].$ Applying Fatou Lemma we obtain that:
        	\begin{equation*}
    	    	0\ge \int_0^T\int_0^T \liminf_\ve \left(\chi_{[t,T]}(s)\left( f'(q_1^\ve)(s)+g'_{\ve_2}(q_2^\ve-q_1^\ve)(s)\right)\right) ds dt =+ \infty,
        	\end{equation*}
        	which is clearly not possible. 
    	\end{proof}
    \end{prop}
    \subsection{Existence of solution of (1)}
    We are now equipped with all the tools and properties needed to show the main result of this section:
    \begin{thm}[Convergence of $q^{\ve}$]
    	\label{thm:convergence_q_ve}
    	For any $\ve_1, \ve_2$ small enough, there exists a subsequence of $q^\ve$ which converges in the $C^{2}$ topology on any compact subset $[\delta, T]$, with $\delta>0$, to a solution $q$ of:
    	\begin{equation*}
    		\begin{cases}
    			\ddot{q}_1 = f'(q_1)+g'(q_2-q_1)\\
    			\ddot{q}_2 = f'(q_2)-g'(q_2-q_1)
    		\end{cases}
    	\end{equation*}
    	with energy $h$ given by 
    	\[
    	h =-f(q_1)-f(q_2)+g(q_2-q_1)\vert_{t=T}<0,
    	\]
    	satisfying $\dot{q}_2(0) = \dot{q}_1(T)=\dot{q}_2(T)=0$ and $q_1(0)=0$.
    	\begin{proof}
    		Recall that the total energy of a solution $q^\ve$ is given by $	h_\ve = -U(q^\ve(T))<0$ (see Proposition \ref{prop:a_priori_estimates_energy}) where $U(q)$ stands for the total potential energy:
    		\[
    		U(q)=f_{\ve_1}(q_1)+f_{\ve_1}(q_2)-g_{\ve_2}(q_2-q_1).
    		\] 
    		Since $q_1^\ve(T)$ and $q_2^\ve(T)$ are uniformly bounded away from $0$, we can assume that $\ve_1$ is so small that $f_{\ve_1}(q_i^\ve(T)) = f(q_i^\ve(T)).$ Since the total energy and the contribution of $f(q_1^\ve(T))+f(q_2^\ve(T))$ are bounded, so is the contribution of $g_{\ve_2}(q_2^\ve(T)-q_1^\ve(T))$.  This implies that there exists a constant $d$, independent on $\ve_1,\ve_2$ (provided that they are small enough), such that $q_2^\ve(T)-q_1^\ve(T)\ge d$. In particular, we have 
    		\begin{equation*}
    			h_\ve = g(q_2^\ve(T)-q_1^\ve(T))- f(q_1^\ve(T))-f(q_2^\ve(T)) <0.
    		\end{equation*}
    		By Proposition  \ref{prop:a_priori_estimates_energy}, $q^\ve$ admits a weakly convergent subsequence and thus converging uniformly and in $L^2$ to some limit function $q = (q_1,q_2)$. Since we have established that $q_2^\ve-q_1^\ve$, $q_1^\ve$ and $q_2^\ve$ are bounded by some positive constants, the right-hand side of \eqref{eq:motion_ve_critical_points} converges uniformly to the bounded function $f'(q_2)-g'(q_2-q_1)$. This fact implies uniform convergence of $\dot q_2^\ve$ too (by Ascoli-Arzelà Theorem). 
    		
    		Differentiating \eqref{eq:motion_ve_critical_points} we obtain:
    		\begin{equation*}
    			\dddot{q}_2^\ve =  f_{\ve_1}''(q_2^\ve)\dot{q}_2^\ve-g''_{\ve_2}(q_2^\ve-q_1^\ve)(\dot q_2^\ve-\dot q_1^\ve).
    		\end{equation*}
    		Again, $f''_{\ve_1}(q_2^\ve)$ and $g''_{\ve_2}(q_2^\ve-q_1^\ve)$ converge uniformly. By Proposition \ref{prop:a_priori_estimates_energy}, $\dot{q}^\ve_1,\dot{q}_2^\ve$ are bounded in $L^2$ and so we obtain that the $\ddot{q}_2^\ve$ are equi-continuous and so they converge uniformly. Thus the limit $q_2$ is $C^2$ and a classical solution of \eqref{eq:motion_ve_critical_points}. 
    		
    		Let us now consider the convergence of $q_1^\ve$. As already mentioned, we have a uniform limit $q_1$. We have to show that $\dot q_1^\ve$ converges uniformly on compact sets of the form $[\delta, T]$, for $\delta>0$ small.  To this aim, let us show that, for any $\delta>0$ there exists $C>0$ such that, for any $\ve_1,\ve_2$ sufficiently small, $q_1^\ve|_{[\delta,T]}\ge C$. If that were not the case, there would be a subsequence $q_1^{\ve_n}$ such that $ q_1^{\ve_n}(\delta)\to 0$. In particular, there would exist $(\omega_n)\sset[0,\delta]$ such that
    		\[
    		\frac{q_1^{\ve_n}(\delta)-q_1^{\ve_n}(0)}{\delta}=\dot{q}_1^{\ve_n}(\omega_n)>\dot{q}_1^{\ve_n}(s)\ge 0,
    		\]
    		for any $s>\delta$. Therefore, $\dot{q}_1^{\ve_n}$ would converge uniformly to 0 on $[\delta,T]$ and, by Poincaré inequality, $q_1^{\ve_n}(T)\to0$ as well, contradicting Proposition \ref{prop:total_collision}.    		
    		
    		The claim we just proved implies that the right-hand side of \eqref{eq:motion_ve_critical_points} converges uniformly to $f'(q_1)+g'(q_2-q_1)$. Applying again Ascoli-Arzelà Theorem we see that $q_1$ is $C^2$ on any $[\delta,T]$ and $q_1^\ve$ converges to $q_1$ in the $C^2$ topology. Thus $q_1$ is a classical solution of \eqref{eq:motion_ve_critical_points} on $[\delta,T]$, for any $\delta>0.$
    	\end{proof}
    \end{thm}
        
        We can specify the result above to the case of homogeneous potentials. Let us take $f(s) = 1/s$ and $g(s) = \mu/s$ with $\mu \in (0,1)$. We obtain:
        \begin{cor}\label{cor:energy}
        	For any negative value of the energy $h$ there exists a \emph{frozen planet orbit}, i.e., a solution of:
        	\begin{equation*}
        		\begin{cases}
        		\ddot{q}_1 = -\frac{1}{q_1^2}-\frac{\mu}{(q_2-q_1)^2} \\
        		\ddot{q}_2 = -\frac{1}{q_2^2}+\frac{\mu}{(q_2-q_1)^2}       		
        	\end{cases}
        	\end{equation*}
           satisfying $\dot{q}_1(T) =q_1(0)= 0$ and $\dot{q}_2(T)=\dot{q}_2(0)=0$.
           \begin{proof}
           	We only need to show that this result holds for any negative value of the energy. Looking at the statement of Proposition \ref{prop:a_priori_estimates_energy} and observing that for these choices of $f$ and $g$ assumption \eqref{hyp:homogeneity} is satisfied with $\al=1$, we get $h_\ve<-c_{\ve_1,\ve_2}/T$. Thus also the energy of the limit obtained in Theorem \ref{thm:convergence_q_ve} is negative. 
           	
           	Since now the potentials are homogeneous, scaling a solution as  $\lambda^{-1}q(\lambda^{3/2} t)$ rescales the energy by a factor $\lambda$, yielding solutions for any negative value of the energy.
           \end{proof}
        \end{cor}

    \section{Zero charge case}
    \label{sec:perturbative}
    Take $\mu\in[0,1] $ and consider the ODE system:
    \begin{equation}
    	\label{eq:ODE_mu}
    	\begin{cases}
    		\ddot{q}_1 = f'(q_1)+\mu g'(q_2-q_1)\\
            \ddot{q}_2 = f'(q_2)-\mu g'(q_2-q_1)
    	\end{cases}.
    \end{equation}
    Theorem \ref{thm:convergence_q_ve} implies that, for any $\mu>0$ there exists a solution $q^\mu = (q_1^\mu,q_2^\mu)$ having $q_1^\mu(0)=\dot{q}^\mu_1(T) =0$ and $\dot{q}^\mu_2(0)=\dot{q}^\mu_2(T)=0$. For $\mu=0$, the system decouples and reduces to two independent $f-$Kepler problems:
    \begin{equation*}
    	\begin{cases}
    		\ddot{q}_1 = f'(q_1)\\
    		\ddot{q}_2 = f'(q_2)
    	\end{cases}.
    \end{equation*}
    Let $\hat{q}$ be a brake orbit for the $f-$Kepler problem on $[0, 2T]$ having $\hat{q}(0)=0$. Define a curve $q =(q_1,q_2)$ in $\mathbb{R}^2$ as:
    \begin{equation}
    	\label{eq:limit_mu_zero}
    	q(t) = (q_1(t),q_2(t)) = (\hat{q}(t),\hat{q}(2T-t)), \quad t \in [0,T].
    \end{equation}
    By construction, we have $\dot{q}_2(0) = q_1(0) =0$ and $q_1(T) = q_2(T)$. The next Proposition shows that, as the charge parameter $\mu$ tends to 0, there exists a sequence of solutions $q^\mu$ converging strongly to $\hat{q}$.
    
    \begin{thm}[Convergence to segments of brake]
    There exists a subsequence of $q^\mu$ which converges uniformly in the $C^2$ topology to the curve $q$ given in \eqref{eq:limit_mu_zero} on any interval $[\delta,T-\delta]$, $\delta>0$. Moreover, there exists a constant $C>0$ such that:
    \[
        q_2^\mu(T)-q_1^\mu(T) \ge C \sqrt[\alpha]{\mu}.
    \]
    	\begin{proof}
    		
    		First let us observe that $q^\mu_2(T)-q^\mu_1(T)\to 0$ as $\mu \to 0$. Indeed, if that were not the case, we could find a subsequence of $(q^\mu _2)$ which would converge in the $C^2$ topology to a solution of $\ddot{x} = f'(x)$, having $\dot{x}(0) = \dot{x}(T) =0$, which does not exist since $f$ is strictly monotone decreasing. Indeed, if $f'(x_0)=0$ for some $x_0$, then $f$ is constant on the whole half-line and thus be zero, thanks to \eqref{hyp:decaying}-\eqref{hyp:monotonicity_convexity}. However, as we have already observed in the proof of Lemma \ref{lem:mountain_pass_geom}, $f(s)\geq f(\bar{s})\left(\frac{\bar{s}}{s}\right)^\al>0$, where $\bar{s}$ is introduced in \eqref{hyp:mountain_pass}. 
    		
    		Next, let us show that $\Vert \dot{q}^\mu \Vert_2^2$ is bounded in $L^2$ and consequently  $q^\mu $ in $H^1$, since $q_1^\mu (0) = 0$ and $q_2^\mu (T)-q_1^\mu (T)\to 0$. Indeed, the value of the action on each $q^\mu $ is uniformly bounded in $\mu$ by some constant $c$ due to the mountain pass structure (see the proof of Proposition \ref{prop:a_priori_estimates_energy} and of Proposition \ref{prop:palais_smale_condition}). Therefore:
    		\begin{equation}
    			\label{eq:proof_critical_point_eq}
    			0 = d_{q^\mu}\mathcal{A}_\mu(q^\mu) \ge \Vert \dot{q^\mu }\Vert^2_2-\al\int_0^T U_{\mu}(q^\mu)\ge \left(1+\frac{\al}{2}\right)\Vert\dot{q}^\mu\Vert_2^2 - \al c,
    		\end{equation}
    	    and so $(\dot{q}^\mu)$ is uniformly bounded in $L^2$. Up to a subsequence, we can assume that the sequence of solutions $(q^\mu)$ converges uniformly and in $L^2$ on $[0,T]$ to a function $q$.
    		
    		Now we prove the $C^2$ convergence on intervals of the form $[\delta,T-\delta]$. Let us observe that $q_1^\mu$ is uniformly strictly concave near $T$. Moreover, thanks to point $iv)$ of Lemma \ref{lemma:convexity_ve_solutions}, $q_2^\mu$ is monotone decreasing and so \[q_1^\mu (t)<q_1^\mu(T)\le q_2^\mu(T)\le q_2^\mu(t). \]
    		Thus the sequence of functions $(q^\mu_2-q_1^\mu)$ cannot converge uniformly to zero on any subinterval of $[0,T].$ A slight modification of the argument given in Theorem \ref{thm:convergence_q_ve} or in Proposition \ref{prop:total_collision} shows that $q_1^\mu$ does not converge uniformly to $0$ on $[0,\delta]$ either.  Applying the Ascoli-Arzelà theorem on the intervals $[\delta,T-\delta]$, we obtain convergence in the $C^{2}$ topology to two solutions $q = (q_1,q_2)$ of $\ddot{x} = f'(x).$     
    		    		
    		We have yet to show that $q_i$ are brake orbits. Let us observe that $q_1^\mu+q_2^\mu$ solves the equation:
    		\begin{equation*}
    			\ddot{q}_1^\mu+\ddot{q}_2^\mu = f'(q_1^\mu)+f'(q_2^\mu).
    		\end{equation*}
    	    The right-hand side is bounded on intervals $[\delta,T]$. This means that we can assume that $q_1^\mu+q_2^\mu$ converges to $q_1+q_2$ in the $C^2$ topology on $[\delta,T]$. Thus $\dot{q}_1(T) = -\dot{q}_2(T)$. Finally, the same argument implies that $q_2^\mu$ converges uniformly in the $C^2$ topology on $[0,T-\delta]$ and so $\dot{q}_2(0)=0$.

    	   Thus, the curve $\hat q$ defined as:
    		\begin{equation*}
    			\hat{q}(t) = \begin{cases}
    				q_1(t) &\text{ if } t\le T\\
    				q_2(2T-t)  &\text{ if } t\ge T
    			\end{cases}
    		\end{equation*} 
    	    determines a $C^1$ trajectory on $[0,2T]$ having $\dot{\hat{q}}(2T)=0$.
    		It follows that $\hat{q}$ is the brake orbit of period $2T$.    		
    		
    		It remains to show that there exists $C>0$ such that $q_2^\mu(T)-q_1^\mu(T) \ge C\mu^{1/\alpha}$. Since we have established convergence of the energies $(h^\mu)$ of the sequence $(q^\mu)$ to the energy  $h$ of $q$, we have:
    		\[
    		\mu g(q_2^\mu(T)-q_1^\mu(T)) \to \dot{q}_1^2(T)>0.
    		\]
    		Moreover, thanks to \eqref{hyp:homogeneity}, the function $g$ satisfies the inequality $g(s)s^\alpha\ge g(1)$ for $s\le 1$. Thus, there exists a constant $C>0$ such that
    		\[q_2^\mu(T)-q_1^\mu(T)\ge C \sqrt[\alpha]{\mu}.\]
    	\end{proof}
    \end{thm}

    \appendix
    \appendixheaderon
    \section{Brake orbits for the $f_{\ve_1}-$Kepler problem}
    
    In this section we briefly discuss the properties of the solution of:
    \begin{equation}
    	\label{eq:brake_kepler}
    	\begin{cases}
    		\ddot{q} = f_{\ve_1}'(q),\\
    		q(0) =0, \, \dot{q}(T)=0,
    	\end{cases}
    \end{equation}
    where $f_{\ve_1}$ has been defined in \eqref{eq:def_f_1} and approximates the attractive potential $f$. Let us consider the space $\mathcal{V} = \{q \in H^1([0,T];\R): q(0)=0\}$, the family of functionals $\mathcal{F}_{\ve_1}$ and $\mathcal{F}$ defined as 
    \begin{equation}
    	\label{eq:smoothed_kepler}
    	\mathcal{F}_{\ve_1} (q)= \int_0^T\frac12 \dot{q}^2+f_{\ve_1}( q), \quad \mathcal{F}(q ) = \int_0^T\frac12 \dot{q}^2+f(\vert q \vert).
    \end{equation}
    Minimisers of $\mathcal{F}$ and $\mathcal{F}_{\ve_1}$ on $\mathcal{V}$ are called \emph{brake orbits} and satisfy \eqref{eq:brake_kepler}.	
    
    \begin{lem}\label{lem:brake}
    	The following assertions hold true:
    	\begin{enumerate}[label = \roman*)]
    		\item $\mathcal{F}$ admits a unique minimiser $\bar q$ in $\mathcal{V}$ which is of class $C^2((0,T])$. This is the unique half brake orbit with minimal period $2T$.
    		\item For any $\ve_1>0$, $\mathcal{F}_{\ve_1}$ has a unique minimiser $q_{\ve_1}$ in $\mathcal{V}$ which is of class $C^2([0,T])$. This is the unique half brake orbit with minimal period $2T$.
    		\item The family $\{q_{\ve_1}\}$ is bounded in $H^1$ and converges strongly to $\bar{q}$.
    		\item Denote by $a_0 = \mathcal{F}(\bar{q})$ and by $a_{\ve_1} =\mathcal{F}_{\ve_1}(q_{\ve_1}).$ Then, $a_{\ve_1}\to a_0$.
    		\item 	Each $\mathcal{F}_{\ve_1}$ satisfies the Palais-Smale condition at any level $c>0.$
    	\end{enumerate}    	    
    	\begin{proof}
    	We first prove points $i)$ and $ii)$. From Lemma \ref{lemma:properties_f_g}, the functions $f_{\ve_1}$ are $C^{1,1}$ on $\mathbb{R}$, the functionals $\mathcal{F}_{\ve_1}$ are $C^1$ and so their critical points are $C^2$. Moreover, the functionals are coercive. Indeed, since $f_{\ve_1}$ is positive, we have:
    		\begin{equation*}
    			\mathcal{F}_{\ve_1}(q) \ge \frac12\int_0^1\vert\dot q\vert^2 dt,
    		\end{equation*}
    		and $\Vert q\Vert_{H^1}\to +\infty $ if and only if $\Vert \dot q\Vert_{2}\to +\infty$ by Poincaré inequality (recall that $q(0)=0$ in $\mathcal{V}$). The functional $\mathcal{F}_{\ve_1}$ is also weakly lower semi-continuous and so minimisers exist for any $\ve_1>0$ by direct methods. Notice that, since we have fixed the starting point and the endpoint is free, any critical point $q_{\ve_1}$ must satisfy $\dot{q}_{\ve_1} (T)=0$. The same argument shows that there exists at least a minimiser for $\mathcal{F}$.
    		
    		For the uniqueness part, let us observe that $T$ can be written in terms of the final position $w =q_{\ve_1}(T).$ Indeed, integrating the energy equation we get:
    	    \begin{align*}
    	     T(w) &= \frac{1}{\sqrt{2}}\int_0^{T(w)} \frac{\dot{q}_{\ve_1}dt}{\sqrt{f_{\ve_1}(q_{\ve_1})-f_{\ve_1}(w)}} =\frac{1}{\sqrt{2}}\int_0^w \frac{dq}{\sqrt{f_{\ve_1}(q)-f_{\ve_1}(w)}}\\
    	     &= \frac{1}{\sqrt{2}}\int_0^1 \frac{ wdq}{\sqrt{f_{\ve_1}(w q)-f_{\ve_1}(w)}}.
    	    \end{align*}
    		Thus, computing the derivative with respect to the final position $w$, we obtain:
    		\begin{equation*}
    			\partial_w T = \frac{1}{2 \sqrt{2}}\int_0^1 \frac{1}{(f_{\ve_1}(wq)-f_{\ve_1}(w))^{3/2}} \left[2 f_{\ve_1}(w q)-2 f_{\ve_1}(w)-q w f'_{\ve_1}(w q)+w f_{\ve_1}'(w)\right].
    		\end{equation*}
    	    So $T$ is strictly monotone in $w$ if the function $\psi(s) = 2 f_{\ve_1}(s)-s f'_{\ve_1}(s)$ has the same property. Computing its derivative we obtain $\psi'(s) = f'_{\ve_1}(s)-s f''_{\ve_1}(s)$, which is always negative in our case.
    	     
    	    Let us prove $iii)$. 
    		By construction, for any curve $q\in\mathcal{V}$, $\mathcal{F}_{\ve_1}(q) \le \mathcal{F}(q)$. It follows that:
    		\begin{equation*}
    			a_{\ve_1} =\min_{q\in \mathcal{V}} \mathcal{F}_{\ve_1}(q) \le \min_{q\in \mathcal{V}} \mathcal{F}(q) =a_0.
    		\end{equation*}
    		Moreover, since $f_{\ve_1}$ is positive, we have:
    		\begin{equation*}
    			\frac12 \Vert \dot{q}_{\ve_1}\Vert^2_2 \le \frac12 \Vert \dot{q}_{\ve_1}\Vert^2_2+	\int_0^T f_{\ve_1}(q_{\ve_1}) \le a_0.
    		\end{equation*}
    		Thus, minimisers of $\mathcal{F}_\ve$ form a bounded subset of $H^1$ and so weakly pre-compact.
    		It remains to prove that $(q_{\ve_1})$ converges strongly in $H^1$ to $\bar{q}$, which is the unique minimiser of $\mathcal{F}$. This is a consequence of Fatou Lemma and uniform convergence:
    		\[
    		\frac12 \Vert \dot{\bar{q}} \Vert^2_2+\int_0^T f(\bar{q})=a_0\ge \liminf_{\ve_1} \left(  \Vert \dot{q}_{\ve_1} \Vert^2_2+\int_0^T f_{\ve_1}(q_{\ve_1})\right)\ge \liminf_{\ve_1}  \Vert \dot{q}_{\ve_1} \Vert^2_2+\int_0^T f(\bar{q}).
    		\]
    		This implies that $	\frac12 \Vert \dot{\bar{q}} \Vert^2_2  = \liminf_{\ve_1}  \Vert \dot{q}_{\ve_1} \Vert^2_2$, proving strong convergence in $H^1$. 
    		
    		Let us show $iv)$. Denoting by $h_{\ve_1}$ the energy of $q_{\ve_1}$, it is not difficult to see that $a_\ve = \Vert \dot{q}_{\ve_1}\Vert_2^2-h_{\ve_1}T$. By the same argument given in Proposition \ref{prop:total_collision} and uniform convergence, we have that $h_{\ve_1}\to \bar{h}$ and so $a_{\ve_1} \to a_0$.
    		
    		It remains to show $v)$. Let $u_n$ be a (PS) sequence for $\mathcal{F}_{\ve_1}$ at a positive level $c>0$. The critical point equation implies:
    		\begin{align*}
    			d_{u_n}\mathcal{F}_{\ve_1}(u_n) &=  \Vert \dot{u}_n\Vert_2^2+ \int_0^Tu_nf_{\ve}'(u_n) \ge \Vert \dot{u}_n\Vert_2^2-\alpha \int_0^Tf_{\ve}(u_n)  \\ &\ge \left(1+\frac{\alpha}{2}\right) \Vert \dot{u}_n\Vert^2 -\alpha \mathcal{F}_{\ve_1}(u_n),
    		\end{align*}
    		where the first inequality follows from point $iii)$ of Lemma \ref{lemma:properties_f_g}. Since $	d_{u_n}\mathcal{F}_\ve(u_n) /\Vert u_n\Vert_2 \to 0$ and we have a Poincaré inequality, $(u_n)$ is a bounded sequence in $H^1$. Thus, up to subsequence, there exists a weak, uniform and $L^2$ limit $u$. To show that it is actually a strong limit, it is enough to notice that:
    		\[
    		d_{u_n}\mathcal{F}_\ve(u-u_n) = \langle \dot{u}_n,\dot u-\dot{u}_n\rangle  + \int_0^Tu_nf_{\ve_1}'(u_n)(u-u_n).
    		\]
    		By the uniform convergence of $u_n$ and  the continuity of $f_{\ve_1}'$, the integral goes to zero and so $\Vert \dot{u}_n\Vert_2 \to\Vert \dot u\Vert_2 $, yielding strong convergence.
    	\end{proof}
    \end{lem}
    
    \section{A mountain-pass Lemma}
    In this appendix we state and prove an \emph{ad-hoc} version of the mountain-pass Lemma. 
    
    \begin{lem}
    	\label{lemma:abstract_mountain_pass}
    	Let $\mathcal{H}$ be a separable Hilbert space, $\mathcal{D}\subset \mathcal{H}$ an open set having smooth boundary and consider a $C^1$ functional $\mathcal{A}\colon \mathcal{D}\to \mathbb{R}$ having Lipschitz gradient. Assume that $\mathcal{A}$ can be extended to $\bar{\mathcal{D}}$ and $\mathcal
    	{A}(\partial \mathcal{D}) = -\infty$.
    	Let $p,q \in \mathcal{D}$, consider the class of paths
    	\begin{equation*}
    		\Gamma = \left\lbrace\gamma\colon[0,1]\to \mathcal{D}\,\middle\lvert\, \gamma(0) =p, \gamma(1)=q\right\rbrace 
    	\end{equation*}
    	and define the value
    	\[
    	c = \inf_{\gamma\in\Gamma} \max_{t\in[0,1]}\mathcal{A}(\gamma(t)).
    	\]
    	Moreover, let us assume that:
    	\begin{itemize}
    		\item $c>\max\{\mathcal{A}(p),\mathcal{A}(q)\}$,
    		\item $\mathcal{A}$ satisfies the $(PS)$-condition at level $c$. 
    	\end{itemize}	
    	Then, $c$ it is a critical value of $\mathcal{A}$.
    	\begin{proof}
    		By contradiction, let us assume that $c$ is not a critical value. since $\mathcal{A}$ satisfy the (PS) condition, this implies that there is no (PS) sequence at level $c$. This means that there exists $\delta>0$ such that, for any $u_n$ with $\mathcal{A}(u_n) \to c$, we have $\Vert d_{u_n}\mathcal{A}\Vert> \delta$. As a consequence, there exists some $\ve>0$ such that $\Vert d_u\mathcal{A}\Vert>2\ve$  for any $u \in \mathcal{A}^{-1}([c-\ve,c+\ve])$. Let us consider the sets:
    		\begin{equation*}
    			X = \mathcal{A}^{-1}([c-2\ve,c+2\ve]), \quad Y = \mathcal{A}^{-1}([c-\ve,c+\ve]).
    		\end{equation*}
    		They are both closed and $\partial \mathcal{D}\cap X =\emptyset$ since $\mathcal{A}(\partial \mathcal{D}) = -\infty$ by assumption. This implies that the function:
    		\begin{equation*}
    			\psi(u) = \frac{\mathrm{dist}(u,X^c)}{\mathrm{dist}(u,X^c)+\mathrm{dist}(u,Y)},
    		\end{equation*}
    		is Lipschitz, vanishes on $X^c$ and is equal to 1 on $Y$. Define the vector field and the associated ODE:
    		\begin{equation*}
    			V = -\frac{\psi \nabla \mathcal{A}}{\Vert\nabla \mathcal{A}\Vert}, \quad \dot{\eta} = V(\eta).
    		\end{equation*}
    		By construction, $V$ is bounded and locally Lipschitz. Thus there exists a well defined continuous flow $\Phi$. It leaves $\mathcal{D}$ invariant since $V\equiv0$ outside $X$ and $\mathcal{A}$ is decreasing on solutions since:
    		\begin{equation*}
    			\frac{d}{dt} \mathcal{A}(\eta(t)) =- \langle \nabla \mathcal{A},\dot{\eta}\rangle = -\psi \Vert\nabla \mathcal{A}\Vert \le 0.
    		\end{equation*}
    		Take $\gamma_n$ a minimising sequence approaching $c$. For $\ve$ sufficiently small, $p,q \in X^c$ and $\Phi_t(\gamma_n)$ still belongs to $\Gamma$. Let $u_n$ be a point in which the maximum is realized. For $n$ sufficiently large, $u_n \in Y$. By definition of $c$ and for $n$ sufficiently large, $c+\ve>\mathcal{A}(\Phi(u_n))\ge c$, however in this case:
    		\begin{equation*}
    			\mathcal{A}(\Phi(u_n))-\mathcal{A}(u_n) = -\int_0^1\Vert\nabla \mathcal{A}(\Phi_s(u_n))\Vert \le -2\ve. 
    		\end{equation*}
    		This implies that $\mathcal{A}(\Phi(u_n)) < c-\ve$, a contradiction.
    	\end{proof}
    \end{lem}
	\bibliography{references}

\begin{thebibliography}{10}

\bibitem{ACZ1987}
A.~Ambrosetti and V.~Coti~Zelati.
\newblock Solutions with minimal period for {H}amiltonian systems in a
  potential well.
\newblock {\em Ann. Inst. H. Poincar\'{e} Anal. Non Lin\'{e}aire},
  4(3):275--296, 1987.

\bibitem{ACZ1990}
A.~Ambrosetti and V.~Coti~Zelati.
\newblock Closed orbits of fixed energy for singular {H}amiltonian systems.
\newblock {\em Arch. Rational Mech. Anal.}, 112(4):339--362, 1990.

\bibitem{AmbCotZel1993}
A.~Ambrosetti and V.~Coti~Zelati.
\newblock {\em Periodic solutions of singular {L}agrangian systems}.
\newblock Progress in Nonlinear Differential Equations and their Applications,
  10. Birkh\"auser Boston Inc., Boston, MA, 1993.

\bibitem{bolotin_degenerate_billiards}
S.~V. Bolotin.
\newblock Degenerate billiards in celestial mechanics.
\newblock {\em Regul. Chaotic Dyn.}, 22(1):27--53, 2017.

\bibitem{Cieliebak_Langmuir}
K.~Cieliebak, U.~Frauenfelder, and M.~Schwingenheuer.
\newblock On {L}angmuir's periodic orbit.
\newblock {\em Arch. Math. (Basel)}, 118(4):413--425, 2022.

\bibitem{Cieliebak_Non_deg}
K.~Cieliebak, U.~Frauenfelder, and E.~Volkov.
\newblock Nondegeneracy and integral count of frozen planet orbits in helium.
\newblock {\em Tunis. J. Math.}, 5(4):713--770, 2023.

\bibitem{Cieliebak_variational}
K.~Cieliebak, U.~Frauenfelder, and E.~Volkov.
\newblock A variational approach to frozen planet orbits in helium.
\newblock {\em Ann. Inst. H. Poincar\'{e} C Anal. Non Lin\'{e}aire},
  40(2):379--455, 2023.

\bibitem{zbMATH06723100}
J.~F{\'e}joz, A.~Knauf, and R.~Montgomery.
\newblock Lagrangian relations and linear point billiards.
\newblock {\em Nonlinearity}, 30(4):1326--1355, 2017.

\bibitem{MR1488438}
J.-M. Rost and G.~Tanner.
\newblock Two-electron atoms: from resonances to fragmentation.
\newblock In {\em Classical, semiclassical and quantum dynamics in atoms},
  volume 485 of {\em Lecture Notes in Phys.}, pages 274--303. Springer, Berlin,
  1997.

\bibitem{Lei_Zhao_shooting}
L.~Zhao.
\newblock Shooting for collinear periodic orbits in the helium model.
\newblock {\em Z. Angew. Math. Phys.}, 74(6):Paper No. 227, 12, 2023.

\end{thebibliography}
	\bibliographystyle{plain}
	
\end{document}